\documentclass[12pt]{amsart}
\usepackage{fullpage,stmaryrd,url,amssymb}
\usepackage[all]{xy}

\newtheorem{theorem}{Theorem}
\newtheorem{question}[theorem]{Question}

\newcommand{\ZZ}{\mathbb{Z}}
\newcommand{\frakm}{\mathfrak{m}}

\DeclareMathOperator{\Aut}{Aut}
\DeclareMathOperator{\cts}{cts}
\DeclareMathOperator{\End}{End}
\DeclareMathOperator{\id}{id}
\DeclareMathOperator{\Sub}{Sub}

\begin{document}

\title{Automorphisms of perfect power series rings}
\author{Kiran S. Kedlaya}
\date{June 29, 2018}
\thanks{Supported by NSF (grant DMS-1501214), UC San Diego (Warschawski Professorship), Guggenheim Fellowship (fall 2015). Some of this work was carried out at MSRI during the fall 2014 research program ``New geometric methods in number theory and automorphic forms'' supported by NSF grant DMS-0932078.
Thanks to Jared Weinstein for formulating the question that led to Theorem~\ref{T:main},
and to Denis Osipov for additional feedback.
}

\begin{abstract}
Let $R$ be a perfect ring of characteristic $p$. We show that the group of continuous $R$-linear automorphisms of the perfect power series ring over $R$ is generated by the automorphisms of the ordinary power series ring together with Frobenius; this answers a question of Jared Weinstein. 
\end{abstract}

\maketitle

\section{Introduction}

Let $R$ be a ring. Let $S = R \llbracket t \rrbracket$ be the ring of formal power series with coefficients in $R$, equipped with the $t$-adic topology, and let $\frakm$ be the ideal of $S$ consisting of series with zero constant term. For each $y \in \frakm$, the formula
\begin{equation} \label{eq:substitution}
\sum_i x_i t^i \mapsto \sum_i x_i y^i
\end{equation}
defines an element $\Sub(y)$ of the algebra $\End^{\cts}_R(S)$ of continuous $R$-linear endomorphisms of $S$ preserving $\frakm$. The resulting map
$\Sub: \frakm \to \End^{\cts}_R(S)$ is inverse to the map $\End^{\cts}_R(S) \to \frakm$ sending $f: S \to S$ to $f(t)$. 
It is well-known (and easy to check) that $\Sub(y)$ is invertible if and only if the coefficient of $t$ in $y$ is a unit in $R$, i.e., 
$y \in R^\times t + \frakm^2$; that is, $\Sub$ identifies $R^\times t + \frakm^2$ with the group $\Aut^{\cts}_R(S)$ of continuous $R$-linear automorphisms of $S$ preserving $\frakm$ (this last condition is automatic if $R$ is reduced).

From now on, assume that $R$ is of characteristic $p>0$ and is perfect, i.e., the Frobenius endomorphism $x \mapsto x^p$ on $R$ is bijective. The analogue of the power series construction in the category of perfect rings is the $t$-adic completion of
$R[t^{1/p}, t^{1/p^2}, \dots]$, which we call $S'$.
The elements of $S'$ may be viewed as formal sums
 $\sum_{i \in \ZZ[p^{-1}]_{\geq 0}} x_i t^i$ with $x_i \in R$ whose support (i.e., the set of $i$ for which $x_i \neq 0$) is either finite or an unbounded increasing sequence.
 
Let $\frakm'$ be the ideal of $S'$ consisting of series with zero constant coefficients.
Then the formula \eqref{eq:substitution} again defines a substitution homomorphism
$\Sub': \frakm' \to \End^{\cts}_R(S')$ which is inverse to the map $\End^{\cts}_R(S') \to \frakm'$ given by evaluation at $t$. In particular, we may construct a commutative diagram
\[
\xymatrix{
\frakm \ar^>>>>{\Sub}@{=}[r] \ar[d] & \End^{\cts}_R(S) \ar[d] \\
\frakm' \ar^>>>>{\Sub'}@{=}[r] & \End^{\cts}_R(S')
}
\]
in which the left vertical arrow is the obvious inclusion. In particular, we get an injective homomorphism $\Aut^{\cts}_R(S) \to \Aut^{\cts}_R(S')$ of groups of continuous $R$-linear automorphisms.

We show that while $\End^{\cts}_R(S')$ is much bigger than $\End^{\cts}_R(S)$, the map of automorphism groups is quite close to being an isomorphism.

\begin{theorem} \label{T:main}
The map $\Aut^{\cts}_R(S) \rtimes \ZZ \to \Aut^{\cts}_R(S')$ taking $n \in \ZZ$ to the map $t \mapsto t^{p^n}$ is an isomorphism of groups.
\end{theorem}

That is, an endomorphism of $S'$ over $R$ is an isomorphism if and only if it carries the uncompleted perfect closure of $S$ into itself. For a related result with perfect closures replaced by algebraic closures, see \cite{kedlaya-temkin}.

Theorem~\ref{T:main} answers a question posed to us by Jared Weinstein, motivated by the following considerations.
Let $S''$ be the $(t_1, t_2)$-adic completion of $R[t_1^{1/p^\infty}, t_2^{1/p^\infty}]$.
By a \emph{(one-dimensional commutative) perfect formal group law} over $R$, we will mean an element
$f \in S''$ satisfying the following conditions.
\begin{enumerate}
\item[(a)]
We have $f(t_2, t_1) = f(t_1, t_2)$.
\item[(b)]
We have $f(t_1) \equiv t_1 \pmod{t_2^{1/p^\infty}}$.
\item[(c)]
In the $(t_1, t_2, t_3)$-adic completion of $R[t_1^{1/p^\infty}, t_2^{1/p^\infty}, t_3^{1/p^\infty}]$, we have $f(f(t_1, t_2), t_3) = f(t_1, f(t_2, t_3))$.
\end{enumerate}
For example, any ordinary (one-dimensional commutative) formal group law over $R$, as an element of $R \llbracket t_1, t_2 \rrbracket$, is a perfect formal group law.

Recall that for every formal group law over a ring of characteristic $p$, the formal multiplication by integers interpolates continuously to a formal action of $\ZZ_p$.
The same holds for perfect formal group laws, and
Theorem~\ref{T:main} implies that for $m \in \ZZ_p^\times$, the formal multiplication-by-$m$ map,
which \emph{a priori} is a perfect power series in one variable, is in fact always an ordinary power series. This suggests a possible affirmative answer to the following question.

\begin{question}
Is every perfect formal group law an ordinary formal group law? That is, is any perfect formal group law contained in $R \llbracket t_1, t_2 \rrbracket$?
\end{question}

It may be possible to gain additional insight into Question~2 by classifying continuous $R$-linear automorphisms of $S''$; however, this approach is complicated by the fact that the map
\[
\Aut^{\cts}_R( R \llbracket t_1, t_2 \rrbracket) \rtimes \ZZ^2 \to \Aut^{\cts}_R(S''),
\]
in which $(n_1, n_2) \in \ZZ^2$ maps to the substitution $t_1 \mapsto t_1^{p^{n_1}}, t_2 \mapsto t_2^{p^{n_2}}$, is far from being surjective. For example, for any $f \in S'$, 
the substitution
\[
t_1 \mapsto t_1, t_2 \mapsto t_2 + f(t_1)
\]
is an automorphism of $S''$ with inverse
\[
t_1 \mapsto t_1, t_2 \mapsto t_2 - f(t_1).
\]


\section{Proof of Theorem~\ref{T:main}}

The remainder of this document consists of the proof of Theorem~\ref{T:main}. The argument is loosely inspired by an analogous calculation of automorphism groups of certain rings of Hahn-Mal'cev-Neumann generalized power series \cite[\S 3]{kp}, although the details turn out to be somewhat different.

By evaluating maps at $t$, we see that the map $\Aut^{\cts}_R(S) \rtimes \ZZ \to \Aut^{\cts}_R(S')$ is injective. To check surjectivity,
let $v_t: S' \to \ZZ[p^{-1}]_{\geq 0}$ denote the $t$-adic valuation, and note that
\[
v_t(\Sub'(y)(z)) = v_t(y) v_t(z) \qquad (y, z \in \frakm').
\]
Consequently, any $y \in \frakm'$ for which $\Sub'(y)$ is an automorphism must satisfy
$v_t(y) \in p^\ZZ$; there is thus no harm in assuming that $v_t(y) = 1$.

It suffices to derive a contradiction assuming that there exist $y,z \in \frakm' \setminus \frakm$ such that $v_t(y) = v_t(z) = 1$ and that $\Sub'(z) \circ \Sub'(y) = \id_{S'}$. (Note that this immediately implies that $\Sub'(y) \circ \Sub'(z) = \id_{S'}$ because $\Sub'(z)$ is injective whenever $z \neq 0$.) Write $y = \sum_i y_i t^i, z = \sum_j z_j t^j$. 
Let $v_p$ denote the $p$-adic valuation on $\ZZ[p^{-1}]$. For $n \in \ZZ$, put
\begin{align*}
a_n &= \min\{i-1: v_p(i) \leq -n, y_i \neq 0\}, \\
b_n &= \min\{i-1: v_p(i) \leq -n, z_i \neq 0\},
\end{align*}
in each case interpreting the minimum over an empty set as $+\infty$.
From our hypotheses,
\begin{equation} \label{eq:ab properties}
a_n, b_n = 0 \,\, (n\leq 0); \,\,
0 < a_1, b_1 < +\infty; \,\,
a_n \leq a_{n+1}, b_n \leq b_{n+1}; \,\,
\lim_{n \to \infty} a_n,
\lim_{n \to \infty} b_n = +\infty.
\end{equation}
Consequently, 
\[
c = \min\left\{\frac{p^n}{p^n-1} a_n, \frac{p^n}{p^n-1} b_n: n=1,2,\dots \right\}
\]
exists, is finite and positive, and is achieved by only finitely many indices. Put
\[
c_n = \frac{p^n-1}{p^n} c \qquad (n \in \ZZ).
\]
Then
\[
a_l \geq c_l, \quad b_m \geq c_m \qquad (l,m \in \ZZ);
\]
there exist maximal indices $l,m$ for which equalities occur;
and these maximal indices are nonnegative and not both zero.
Moreover, if $l>0$, then we must have $v_p(a_l) = -l$, as otherwise we would have the contradiction
\[
c_l = a_l = a_{l+1} \geq \frac{p^{l+1}-1}{p^{l+1}} c > \frac{p^l-1}{p^l}c = c_l;
\]
similarly, if $m>0$, then $v_p(b_m) = -l$.
Since we either have $a_l = c_l$ for some $l>0$ or $b_m = c_m$ for some $m>0$,
we may deduce that for all $n>0$, $c_n > 0$ and $v_p(c_n) = -n$.

Since $l+m > 0$, we have $c_{l+m} > 0$, so the coefficient of $t^{1+c_{l+m}}$ in 
$t = (\Sub'(z) \circ \Sub'(y))(t) = \sum_{i>0} y_i z^i$ must be zero.
To obtain the desired contradiction, it will thus suffice to verify that the coefficient of $t^{1+c_{l+m}}$ in $y_i z^i$ is nonzero for exactly \emph{one} value of $i$; we check this by distinguishing options for $d = -v_p(i)$.
\begin{itemize}
\item
For $d \geq l+m$, we have
\[
v_t(y_i z^i) = v_t(y_i t^i) \geq 1 + a_d \geq 1 + c_d \geq 1 + c_{l+m}.
\]
For the coefficient of $t^{1+c_{l+m}}$ in $y_i z^i$ to be nonzero,
this chain of inequalities must become a chain of equalities, yielding
\[
i = 1 + a_d, \qquad d \leq l, \qquad d = l+m.
\]
Since $m \geq 0$, this is only possible if
$d = l$, $m = 0$, $i = 1+a_l$;
in this case, the coefficient of 
$t^{1 + c_{l+m}}$ in $y_i z^i$ is the nonzero value $y_{1+a_l} z_1^{1+a_l}$.
\item
For $d < l+m$, we have
\[
y_i z^i = y_i (z^{p^{-d}})^{ip^d} = y_i z_1^i \left( 1 + \sum_{j>1} (z_j/z_1) t^{p^{-d}(j-1)} \right)^{ip^d} t^{i}.
\]
By the definition of $b_n$, the sum over $j$ can be rewritten as
\[
\text{(element of $S^{p^{-l-m+1}}$)} + 
\text{(nonzero element of $R$)} \cdot t^{p^{-d} b_{l+m-d}} 
+ \text{(higher order terms)}.
\]
Since $v_p(ip^d) = 0$ and $t^i \in S^{p^{-l-m+1}}$, the binomial expansion yields
\[
y_i z^i = \text{(element of $S^{p^{-l-m+1}}$)} + 
\text{(nonzero element of $R$)} \cdot t^{i + p^{-d} b_{l+m-d}} 
+ \text{(higher order terms)}.
\]
Since $v_p(c_{l+m}) = -l-m$, the coefficient of $t^{1+c_{l+m}}$ in any element of $S^{p^{-l-m+1}}$ is zero. On the other hand, we have
\[
i+ p^{-d} b_{l+m-d} \geq 1 + a_d + p^{-d} b_{l+m-d} \geq 1 + c_d + p^{-d} c_{l+m-d} = 1 + c_{l+m}.
\]
For the coefficient of $t^{1+c_{l+m}}$ in $y_i z^i$ to be nonzero,
this chain of inequalities must become a chain of equalities, yielding
\[
i = 1 + a_d, \qquad d \leq l, \qquad l+m-d \leq m.
\]
This is only possible if
$d = l$, $i = 1+a_l$, $m>0$;
in this case, the coefficient of 
$t^{1 + c_{l+m}}$ in $y_i z^i$ is the nonzero value $y_{1+a_l} z_1^{1+a_l}$.
\end{itemize}

Since exactly one of the two boundary cases can occur (depending on whether $m=0$ or $m>0$),
this yields the desired contradiction.

\end{document}